\documentclass[12pt]{article}

\usepackage[a4paper, total={6in, 9in}]{geometry}

\usepackage{mathptmx}
\usepackage{times}

\usepackage{amsmath,amsfonts,amssymb,amsthm}
\usepackage{mathrsfs}
\usepackage{bm}

\usepackage{tikz}
\usetikzlibrary{calc, trees, positioning}
\tikzstyle{every node}=[circle,inner sep=1pt,fill=white!60]
\tikzstyle{tn}=[shape=circle, draw, color=black!70]
\tikzstyle{tn1}=[shape=circle, draw]
\tikzstyle{tn2}=[shape=circle, draw,inner sep=1.5pt]
\tikzstyle{tn5}=[shape=circle, inner sep=1.6pt,draw, color=black!70]
\tikzstyle{tn3}=[shape=rectangle, draw,inner sep=1.5pt]
\tikzstyle{tn4}=[shape=rectangle, draw,inner sep=1.5pt,color=black!70]
\tikzstyle{marke}=[shape=circle,minimum size=0.1cm, draw,blue]

\usepackage{verbatim}
\usepackage{soul}
\usepackage{hyperref}
\usepackage[numbers,sort&compress]{natbib}

\newtheorem{thm}{Theorem}[section]

\numberwithin{equation}{section}

\def\qed{\nopagebreak\hfill{\rule{4pt}{7pt}}}

\newcommand\exc{{\mathrm{exc}}}
\newcommand\asc{\mathrm {asc}}

\newcommand\drop{\mathrm{drop}}
\newcommand\des{\mathrm{des}}

\newcommand\cyc{\mathrm{cyc}}
\newcommand\LR{\mathrm{lrmin}}
\newcommand\RL{\mathrm{rlmin}}
\newcommand\xleaf{\mathrm{xleaf}}
\newcommand\yleaf{\mathrm{yleaf}}

\allowdisplaybreaks
\begin{document}

\begin{center}
{\Large\bf 
Increasing Binary Trees and 
 
the $(\alpha,\beta)$-Eulerian  Polynomials
 }

\vskip 6mm

William Y.C. Chen$^1$ and Amy M. Fu$^2$

\vskip 3mm

$^{1}$Center for Applied Mathematics, KL-AAGDM \\
Tianjin University\\
Tianjin 300072, P.R. China

\vskip 3mm

$^{2}$School of Mathematics\\
Shanghai University of Finance and Economics\\
Shanghai 200433, P.R. China

\vskip 3mm

Emails: { $^1$chenyc@tju.edu.cn, $^{2}$fu.mei@mail.shufe.edu.cn}

\vskip 6mm

\end{center}

\begin{abstract}  
In light of the grammar given by Ji for the $(\alpha,\beta)$-Eulerian     polynomials introduced by Carlitz and Scoville, we provide a labeling scheme for increasing binary trees. In this setting, we obtain a combinatorial interpretation 
of the $\gamma$-coefficients of the $\alpha$-Eulerian polynomials in terms of forests of planted 0-1-2-plane trees, 
which specializes to a combinatorial interpretation of the $\gamma$-coefficients of the derangement polynomials in the same spirit. 
By means of a decomposition of an increasing binary tree into a forest,
we find combinatorial interpretations of the sums involving two identities of Ji, one of which can be viewed as $(\alpha,\beta)$-extensions of the formulas of Petersen and Stembridge. 
\end{abstract}



\section{Introduction}

The objective of this paper is to explore a labeling scheme 
for increasing binary trees as an alternative combinatorial
interpretation of the  $(\alpha,\beta)$-Eulerian polynomials
introduced by Carlitz and Scoville  \cite{Carlitz-Scoville-1974}. 
 A grammatical treatment of these polynomials
has been given by Ji \cite{Ji-2023-A} 
via a  labeling scheme for permutations.
Employing the grammatical calculus, Ji obtained 
  $(\alpha, \beta)$-extensions 
of the formulas of Petersen and Stembridge. 

We begin with a combinatorial setting 
of the $(\alpha,\beta)$-Eulerian polynomials
in terms of increasing binary trees. 
Based on an equivalent definition of Ji
relying on the number of left-to-right minima and
the number of right-to-left minima of a permutation,
we observe that two particular leaves 
of an increasing binary tree, called 
the $a$-leaf and the $b$-leaf, play a special role.
Then we move on to define
the $\alpha$-vertices and the $\beta$-vertices, and 
add the $\alpha$-labels and the
$\beta$-labels to certain internal vertices,
while adopting 
the $(x,y)$-labeling for the leaves, 
as given in  \cite{Chen-Fu-2022} for the bivariate Eulerian polynomials.

 In fact, the
two special leaves (the $a$-leaf and the
$b$-leaf) can be considered as two poles to stretch
a binary tree aligned on a horizontal line,
which is reminiscent of
the decomposition of a doubly rooted tree into a linear
order of rooted trees in Joyal's proof of
Cayley's formula \cite{Joyal}. 
More precisely, with these two special vertices
at disposal, an 
increasing binary tree can be decomposed 
into a forest of planted   increasing
binary trees.
Such a decomposition gives rise to
a combinatorial interpretation
of the $\gamma$-coefficients of the $\alpha$-Eulerian polynomials
in terms of forests of planted 0-1-2-plane trees.
An interpretation 
in the permutation setting has been 
given by Ji-Lin \cite{Ji-Lin-2023} by
devising a group action.

The idea of the labeling scheme for the
$(\alpha,\beta)$-Eulerian polynomials can be
adapted to a grammar of Dumont related to the
derangement polynomials. In this setting, we 
are led to 
a combinatorial interpretation of the $\gamma$-coefficients
of the derangement polynomials and the $q$-derangement polynomials (with respect to the 
number of cycles),
in terms of forests of planted increasing
0-1-2-plane trees, where the exponents of $q$ are
connected with the number of components of a forest.
This topic
has been extensively studied, see, for example,
\cite{LZ-2015, LSW-2012, MMYY-2024, MMYZ-2018, SW-2020, SZ-2010, Sun-Wang-2014}. 

 The grammatical labelings of increasing
 binary trees make it possible to give
 combinatorial interpretations of the sums
 involving the identities of Ji.  
 We first realize that
the number of  interior peaks 
of a permutation can be read off from a labeling  of increasing
binary trees.  
For the rest, the decomposition of an increasing binary tree  is the key ingredient all along.

\section{The $(\alpha,\beta)$-Eulerian polynomials}

For $n\geq 1$, let $[n]=\{1,2,\ldots,n\}$.
Given a permutation $\sigma=\sigma_1 \cdots \sigma_n$ of $[n]$, 
an index $i$ ($2\leq i \leq n$) is called an ascent if $\sigma_{i-1}<\sigma_{i}$,
and an index 
$i$ ($1\leq i \leq n-1$) is called a descent if $\sigma_i>\sigma_{i+1}$. Let $\asc(\sigma)$ and $\des(\sigma)$ denote the number of ascents and the number of descents of $\sigma$, respectively.

Carlitz and Scoville \cite{Carlitz-Scoville-1974} introduced an extension of the
bivariate Eulerian polynomials, denoted by 
$A_n(x,y\,|\,\alpha, \beta)$, which are called 
the $(\alpha,\beta)$-Eulerian polynomials by 
Ji \cite{Ji-2023-A}. 
They are defined by
\begin{equation}
A_n(x,y  \,|\,  \alpha,\beta)=
\sum_{\sigma\in{S}_{n+1}}x^{{\mathrm {asc}}(\sigma)}y^{{\des}(\sigma)}\alpha^{{\rm lrmax}(\sigma)-1}\beta^{{\rm rlmax}(\sigma)-1},
\end{equation}
where $S_{n+1}$ is the set of 
permutations of $[n+1]$,
${\rm lrmax}(\sigma)$ and ${\rm rlmax}(\sigma)$ denote 
the number of left-to-right maxima and the number of right-to-left
maxima of $\sigma$, respectively.

By taking complement of a permutation and exchanging
the roles of $x$ and $y$, Ji \cite{Ji-2023-A}
presented an equivalent definition
\begin{equation}\label{Anxy}
A_n(x,y\,|\,\alpha,\beta)=
\sum_{\sigma\in{S}_{n+1}}x^{{\mathrm{des}}(\sigma)}y^{{\asc}(\sigma)} \alpha ^{{\LR}(\sigma)-1}\beta^{{\RL}(\sigma)-1},
\end{equation}
where ${\LR}(\sigma)$ and ${\RL}(\sigma)$ denote 
the number of left-to-right minima and the number of right-to-left
minima of $\sigma$, respectively. 
The initial values of $A_n(x,y\,|\,\alpha,\beta)$
are given below,
\begin{eqnarray*}
  A_0(x,y\,|\,\alpha,\beta) &  =  & 1, \\[3pt]
  A_1(x,y\,|\,\alpha,\beta) &  =  &x \beta + y \alpha , \\[3pt]
  A_2(x,y\,|\,\alpha,\beta) &  =  & xy \alpha + xy\beta
  + 2xy \alpha \beta + x^2 \beta^2 + y^2\alpha^2.
  \end{eqnarray*}

Ji \cite{Ji-2023-A} found a context-free grammar 
  for the $(\alpha,\beta)$-Eulerian polynomials, which
  can be paraphrased as
\begin{equation}
    G=\{ a \rightarrow \alpha a y, \;
         b \rightarrow \beta b x , \;
         x \rightarrow xy, \;
         y \rightarrow xy \}.
\end{equation}
By providing a grammatical labeling for permutations, it has
been shown that the $(\alpha,\beta)$-Eulerian polynomials
can be generated by the above grammar.

\begin{thm}[Ji]
Let $D$ denote the formal derivative of the
  grammar $G$. For $n\geq 0$, we have 
\begin{equation}
   D^n (ab) = ab\, A_{n}(x,y\,|\, \alpha, \beta).
\end{equation}
\end{thm}

 As is well-known, permutations are in one-to-one
correspondence with increasing binary trees, we find that
endowed with a suitable labeling scheme increasing binary trees
are conducive to a combinatorial understanding of 
the $(\alpha,\beta)$-Eulerian polynomials. For this purpose,
we shall introduce the $(a,b,\alpha,\beta)$-labeling as
described below.

\subsection{The $(a,b,\alpha,\beta)$-labeling}

Let $n\geq 1$, and let $T$ be an
increasing binary tree on $[n]$, where
$n\geq 1$. Consider the left child of the root.
If it is a leaf, we call it the leftmost leaf of $T$.
If not, we restrict to the left subtree of $T$ and continue to
seek the leftmost leaf. Eventually, we
end up with  the
leftmost leaf of $T$. The rightmost leaf is 
 defined in the same way. Now
 we label leftmost leaf of $T$  by $a$ and label
the
rightmost leaf of $T$ by $b$.

\begin{figure}[!ht]
\begin{center}
\begin{tikzpicture}[scale=0.8]
\node [tn,label=90:$1$]{}[grow=down]
	[sibling distance=36mm,level distance=16mm]
    child {node [tn,label=180:{$4(\alpha)$}](four){}
       [sibling distance=25mm,level distance=13mm]
    child {node [tn,label=180:{$8(\alpha)$}](eight){}
     [sibling distance=14mm,level distance=13mm]
     child {node [tn1,label=-90:{}](eightl){}}
     child {node [tn1,label=0:{}](eightr){}}
     }
      child {node [tn,label=0:{$6$}](six){}
     [sibling distance=14mm,level distance=13mm]
     child {node [tn,label=180:{$9$}](nine){}
     child {node [tn1, label=-90:{}](ninel){}}
     child {node [tn1, label=-90:{}](niner){}}
     }
     child {node [tn1,label=0:{}](sixr){}}
     }
     }
     child {node [tn,label=0:{$2(\beta)$}](twoa){}
     [sibling distance=22mm,level distance=13mm]
      child {node [tn1,label=-90:{$x$}](twoal){}}
    child {node [tn,label=0:{$3(\beta)$}](threea){}
       [sibling distance=20mm,level distance=13mm]
    child {node [tn,label=180:{$5$}](five){}
       [sibling distance=14mm,level distance=13mm]
      child {node [tn1,label=-90:{}](fivel){}}
    child {node [tn1,label=-90:{}](fiver){}}
    }
    child {node [tn,label=0:{$7(\beta)$}](sevena){}
         [sibling distance=14mm,level distance=13mm]
    child {node [tn1,label=-90:{}](sevenl){}}
    child {node [tn1,label=-90:{$b$}](sevenr){}}
    }
     }};
    \node [below=3pt] at (eightl){$a$}; 
     \node [below=3pt] at (eightr){$y$}; 
     \node [below=3pt] at (sixr){$y$};
     \node [below=3pt] at (ninel){$x$};
     \node [below=3pt] at (niner){$y$};
     \node [below=3pt] at (fivel){$x$}; 
     \node [below=3pt] at (fiver){$y$};
    \node [below=3pt] at (sevenl){$x$};
    \node [below=3pt] at (sevenr){$b$};
\end{tikzpicture}
\end{center}
\caption{An example for the $(a,b,\alpha,\beta)$-labeling.}
\label{aabbf}
\end{figure}
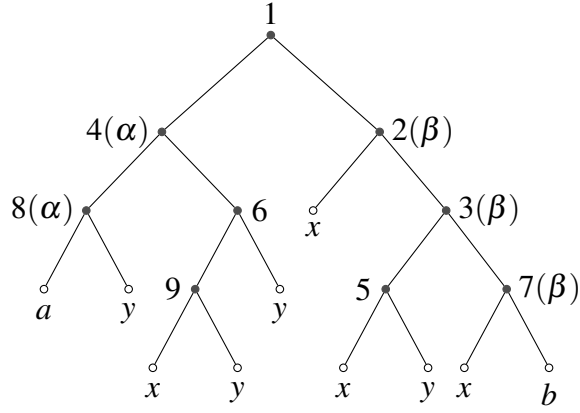

Next, the $\alpha$-vertices and the
$\beta$-vertices are defined as follows.
Each vertex on the path from the root
to the $a$-leaf (other than the root   
and the $a$-leaf) 
is labeled by $\alpha$, which we call an $\alpha$-vertex. 
Each vertex on the path from the root 
to the $b$-leaf (other than the root  and
the $b$-leaf)  is labeled by $\beta$, which 
we call a $\beta$-vertex. The rest of
the leaves are labeled like the usual $(x,y)$-labeling, that is, a left leaf is labeled by 
$x$ and a right leaf is labeled by $y$. 
It can be readily seen that
a pair of sibling leaves labeled by $x$
and $y$ correspond to an interior 
peak of a permutation. For example,  
Figure \ref{aabbf} demonstrates an $(a,b,\alpha,\beta)$-labeling of an increasing binary tree on $[ 9]$,  where the corresponding
permutation reads
\[ 8 \ 4 \ 9 \ 6 \  1 \ 2\ 5\ 3\ 7. \]

The following theorem shows that 
the $(\alpha,\beta)$-Eulerian polynomials 
have a combinatorial interpretation in terms
of increasing binary trees. For an
increasing binary tree $T$, we use $w(T)$ denote
the weight of $T$ with respect to the
$(a,b,\alpha,\beta)$-labeling, that is,
the product of the grammatical labels. 
For instance, the weight of the increasing
binary in Figure \ref{aabbf} equals
$abx^4y^4\alpha^2\beta^3$.

In view of the correspondence between permutations
and increasing binary trees, we see that
the $\alpha$-vertices together with the root
are the
left-to-right minima of the corresponding permutation,
whereas the  $\beta$-vertices together with 
the root
are the right-to-left minima of the 
corresponding permutation.  Then we arrive at the following combinatorial expansion. 
To be consistent with the meanings of $x$ and $y$ in the
$(a,b,\alpha,\beta)$-labeling, we use $A^*_n(x,y\,|\,\alpha,\beta)$
to denote $A_n(y,x\,|\,\alpha,\beta)$.

\begin{thm}
For $n\geq 1$, we have
\begin{equation}
    ab\,A^*_n(x,y\,|\,\alpha,\beta)= \sum_{T} w(T) ,
\end{equation}
where the sum ranges over the set of
increasing binary trees on $[n+1]$
with the $(a,b,\alpha,\beta)$-labeling.
\end{thm}

\subsection{A decomposition}

The $(a,b,\alpha,\beta)$-labeling leads us to 
consider a decomposition of an increasing binary
tree into a forest of planted increasing binary trees, which can be used to divide the set of 
increasing binary trees into classes
relative to the labeling scheme.  By a 
planted increasing binary tree we mean 
an increasing tree structure consisting of a
single root or a root with an increasing binary
tree as a subtree. In the usual sense, 
a planted plane tree is either a single root
or a plane tree for which the root has only one
child.

Next, we introduce
a decomposition of an
increasing binary tree $T$ on $[n]=
\{ 1,2, \ldots, n\}$ 
into a forest of 
planted increasing binary trees
rooted at the $\alpha$-vertices and the $\beta$-vertices. The resulting forest is called the supporting
forest of $T$. For an increasing binary
tree $T$ on $[n]$, its supporting forest
is on the set $[2,n]=\{2, 3, \ldots, n\}$.

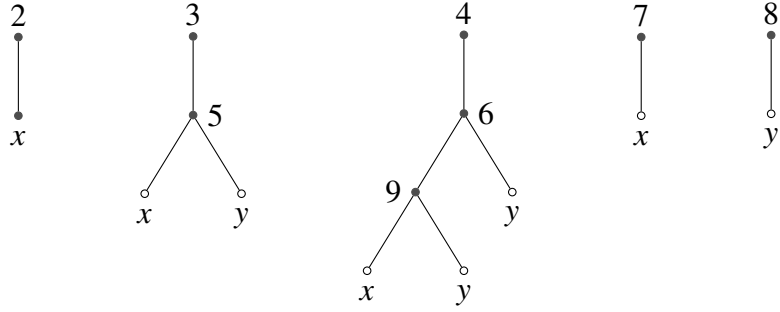
\begin{figure}[!ht]
\begin{center}
\raisebox{60pt}{\begin{tikzpicture}[scale=0.8]
\node [tn,label=90:$2$]{}[grow=down]
[sibling distance=14mm,level distance=13mm]
    child {node [tn,label=-90:{}](two){}};
    \node [below=3pt] at (two){$x$}; 
\end{tikzpicture}}
\hspace{30pt}
\raisebox{29pt}{\begin{tikzpicture}[scale=0.8]
\node [tn,label=90:$3$]{}[grow=down]
[sibling distance=14mm,level distance=13mm]
    child {node [tn,label=0:{$5$}](five){}
    [sibling distance=16mm,level distance=13mm]
    child {node [tn1,label=-90:{}](fivel){}}
    child {node [tn1,label=-90:{}](fiver){}}};
    \node [below=3pt] at (fivel){$x$}; 
    \node [below=3pt] at (fiver){$y$}; 
\end{tikzpicture}}
\hspace{30pt}
\begin{tikzpicture}[scale=0.8]
\node [tn,label=90:$4$]{}[grow=down]
[sibling distance=14mm,level distance=13mm]
    child {node [tn,label=0:{$6$}](six){}
    [sibling distance=16mm,level distance=13mm]
    child {node [tn,label=180:{$9$}](nine){}
    [sibling distance=16mm,level distance=13mm]
    child {node [tn1,label=180:{}](ninel){}}
    child {node [tn1,label=180:{}](niner){}}}
    child {node [tn1,label=-90:{}](sixr){}}};
    \node [below=3pt] at (ninel){$x$}; 
    \node [below=3pt] at (niner){$y$}; 
    \node [below=3pt] at (sixr){$y$}; 
\end{tikzpicture}
\hspace{30pt}
\raisebox{60pt}{\begin{tikzpicture}[scale=0.8]
\node [tn,label=90:$7$]{}[grow=down]
[sibling distance=14mm,level distance=13mm]
    child {node [tn1,label=-90:{}](seven){}};
    \node [below=3pt] at (seven){$x$}; 
\end{tikzpicture}}
\hspace{30pt}
\raisebox{59pt}{\begin{tikzpicture}[scale=0.8]
\node [tn,label=90:$8$]{}[grow=down]
[sibling distance=14mm,level distance=13mm]
    child {node [tn1,label=-90:{}](eight){}};
    \node [below=3pt] at (eight){$y$}; 
\end{tikzpicture}}
\end{center}
\caption{A supporting  forest on
$[2,n]$ with inherited labels.}
\label{aabbforest}
\end{figure}

If we arrange the components of a
supporting forest in the increasing order of their roots,  then Figure \ref{aabbforest} is an 
exhibition of the supporting forest of the increasing tree in Figure \ref{aabbf}.

The structure of a supporting forest can be
used to divide the set of increasing binary
trees on $[2,n]$ into classes whose total weight
can be readily characterized. 
To this end, we define the weight of a 
supporting forest by the following rules. First,
we suppress the leaf of a single root. 
\begin{enumerate}
\item A single root is assigned the weight $x\beta+y\alpha$. 
\item A root with a child has weight $\alpha+\beta$.
\item Any leaf has the weight (or label)
inherited   from 
 the original increasing binary tree. 
\end{enumerate}

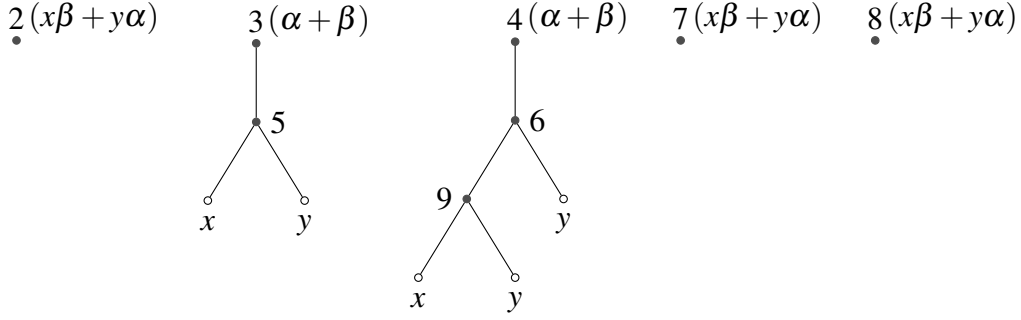
\begin{figure}[!ht]
\begin{center}
\raisebox{103pt}{\begin{tikzpicture}[scale=0.8]
\node [tn,label=90:$2$]{};
\end{tikzpicture}}
\hspace{-5pt}\raisebox{110pt}{${\textstyle (x\beta+y\alpha)}$}
\hspace{7pt}
\raisebox{29pt}{\begin{tikzpicture}[scale=0.8]
\node [tn,label=90:$3$]{}[grow=down]
[sibling distance=14mm,level distance=13mm]
    child {node [tn,label=0:{$5$}](five){}
    [sibling distance=16mm,level distance=13mm]
    child {node [tn1,label=-90:{}](fivel){}}
    child {node [tn1,label=-90:{}](fiver){}}};
    \node [below=3pt] at (fivel){$x$}; 
    \node [below=3pt] at (fiver){$y$}; 
\end{tikzpicture}}
\hspace{-23pt}\raisebox{109pt}{${\textstyle (\alpha+\beta)}$}
\hspace{7pt}
\begin{tikzpicture}[scale=0.8]
\node [tn,label=90:$4$]{}[grow=down]
[sibling distance=14mm,level distance=13mm]
    child {node [tn,label=0:{$6$}](six){}
    [sibling distance=16mm,level distance=13mm]
    child {node [tn,label=180:{$9$}](nine){}
    [sibling distance=16mm,level distance=13mm]
    child {node [tn1,label=180:{}](ninel){}}
    child {node [tn1,label=180:{}](niner){}}}
    child {node [tn1,label=-90:{}](sixr){}}};
    \node [below=3pt] at (ninel){$x$}; 
    \node [below=3pt] at (niner){$y$}; 
    \node [below=3pt] at (sixr){$y$}; 
\end{tikzpicture}
\hspace{-23pt}\raisebox{110pt}{${\textstyle (\alpha+\beta)}$}
\hspace{7pt}
\raisebox{103pt}{\begin{tikzpicture}[scale=0.8]
\node [tn,label=90:$7$]{};
\end{tikzpicture}}
\hspace{-5pt}\raisebox{110pt}{${\textstyle (x\beta+y\alpha)}$}
\hspace{7pt}
\raisebox{103pt}{\begin{tikzpicture}[scale=0.8]
\node [tn,label=90:$8$]{};
\end{tikzpicture}}
\hspace{-5pt}\raisebox{110pt}{${\textstyle (x\beta+y\alpha)}$}
\end{center}
\caption{A supporting  forest with updated labels.}
\label{updated}
\end{figure}

The updated labeling of a supporting forest
is illustrated in Figure \ref{updated}.

Since the root of a component of a supporting forest
can be either an $\alpha$-vertex or a $\beta$-vertex,
we are led to the following expansion,
where the underlying set of the
supporting forests has been rescaled down to
$[n]$. 

\begin{thm}
    For $n\geq 0$, 
        $A^*_n(x,y\,|\,\alpha,\beta)$
    equals the total weight of supporting forests on $[n]$.
\end{thm}

Now we further classify supporting forests via
a group action. We say that two supporting forests
are in the same class if one can be obtained from another by 
swapping a leaf with its non-leaf sibling. Therefore,
such a class of supporting forests can be represented
by a forest of planted 0-1-2 plane trees (without
external leaves), bearing
the following labeling rules:
\begin{enumerate}
\item A single root is endowed with a weight $x\beta+y\alpha$. 
\item A root with a child has weight $\alpha+\beta$.
\item A degree one non-root vertex (a nonroot
vertex with exactly one  child) has weight
$x+y$.
\item A leaf has weight $xy$.
\end{enumerate}

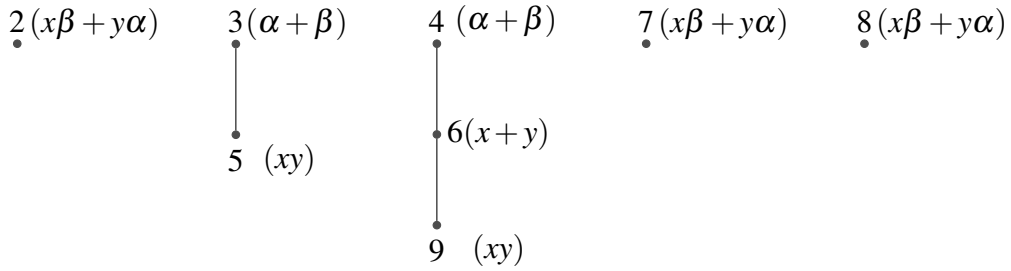
\begin{figure}[!ht]
\begin{center}
\raisebox{83pt}{\begin{tikzpicture}[scale=0.8]
\node [tn,label=90:$2$]{};
\end{tikzpicture}}
\hspace{-5pt}\raisebox{89pt}{${\textstyle (x\beta+y\alpha)}$}
\hspace{16pt}
\raisebox{34pt}{\begin{tikzpicture}[scale=0.8]
\node [tn,label=90:$3$]{}[grow=down]
    child {node [tn,label=0:{}](five){}
    };
    \node [below=3pt] at (five){$5$};  
\end{tikzpicture}}
\hspace{-6pt}\raisebox{89pt}{${\textstyle (\alpha+\beta)}$}
\hspace{-35pt}\raisebox{38pt}{${\textstyle (xy)}$}
\hspace{33pt}
\begin{tikzpicture}[scale=0.8]
\node [tn,label=90:$4$]{}[grow=down]
    child {node [tn,label=0:{$6(x+y)$}](six){}
    child {node [tn,label=-90:{}](nine){}}}; 
    \node [below=3pt] at (nine){$9$}; 
\end{tikzpicture}
\hspace{-41pt}\raisebox{90pt}{${\textstyle (\alpha+\beta)}$}
\hspace{-35pt}\raisebox{4pt}{${\textstyle (xy)}$}
\hspace{33pt}
\raisebox{83pt}{\begin{tikzpicture}[scale=0.8]
\node [tn,label=90:$7$]{};
\end{tikzpicture}}
\hspace{-5pt}\raisebox{89pt}{${\textstyle (x\beta+y\alpha)}$}
\hspace{16pt}
\raisebox{83pt}{\begin{tikzpicture}[scale=0.8]
\node [tn,label=90:$8$]{};
\end{tikzpicture}}
\hspace{-5pt}\raisebox{89pt}{${\textstyle (x\beta+y\alpha)}$}
\end{center}
\caption{A forest of planted 0-1-2-plane trees.}
\label{plane}
\end{figure}

Figure \ref{plane} is an illustration of a
forest of planted 0-1-2-plane trees. 
The above classification 
implies the following expansion of the
$(\alpha,\beta)$-Eulerian polynomials.

\begin{thm}
    For $n\geq 0$, 
     $A^*_n(x,y\,|\,\alpha,\beta)$
    equals
    the total weight of forests of 
    planted 0-1-2-plane trees on $[n]$.
\end{thm}

\subsection{The $\alpha$-Eulerian polynomials }

When $\alpha=\beta$, the $(\alpha,\beta)$-Eulerian polynomials
are called the $\alpha$-Eulerian polynomials in
\cite{Ji-2023-A}, denoted by  
$A_n(x,y\,|\,\alpha)$. Likewise, we use
$A_n^*(x,y \, | \, \alpha)$ to denote
$A_n(y,x \, | \, \alpha)$.
The labeling scheme of the corresponding trees is  called the $(a,b,\alpha)$-labeling. That is, all $\beta$-vertices are  labeled by $\alpha$ as well. 
By a transformation of grammars, it is easy
to see that these polynomials are $\gamma$-positive. Recall that  by a $\gamma$-expansion, we mean an expansion
in $x+y$ and $xy$. A bivariate polynomial
is called $\gamma$-positive if the coefficients
of the $\gamma$-expansion are all nonnegative. 
Evidently, the usual
notion of $\gamma$-positivity for polynomials
in $x$ is equivalent to the bivariate formulation, whereas
we do need both variables $x$ and $y$ as far as
the grammar is concerned.

Ji and Lin \cite{Ji-Lin-2023} provided a combinatorial proof
of the $\gamma$-coefficients by
via a  group action on permutations. 
With the help of the $(a,b,\alpha)$-labeling, we obtain an alternative combinatorial
interpretation of the $\gamma$-coefficients
in terms of forests of planted 
0-1-2-plane trees.

Setting $\alpha=\beta$, the
previous weight assignment reduces to 
the following rules for the $\alpha$-Eulerian
polynomials. For a forest $F$ of planted
increasing 0-1-2-plane trees, we
have the following rules:
\begin{enumerate}
\item 
      A single root has weight $\alpha(x+y)$.
      
\item Other roots have weight
       $2 \alpha$.

    \item If a non-root vertex has only one child, it has weight $x+y$.

    \item A leaf has weight $xy$. 
\end{enumerate}

\begin{thm}For $n\geq 1$, the $\alpha$-Eulerian polynomial
$A^*_n(x,y\,|\,\alpha)$
has the $\gamma$-expansion 
\[ \sum_{F}  w(F),\]
where the sum ranges over forests of 
planted 0-1-2-plane trees on $[n]$. 
\end{thm}

\subsection{The derangement polynomials}

As a special case of the $\gamma$-expansion
of the $\alpha$-Eulerian polynomials, we come
to the $\gamma$-expansion of the derangement polynomials.

Given a permutation $\sigma$ in the cycle notation, 
assume that the minimum element of each cycle appears at
the end, and the cycles are arranged in the 
increasing order of their minimum elements. For example, 
 $(8\, 4\, 9\, 6\, 1) \;(2 )\; (5\, 3 )\; (7) $ 
 is  a permutation of $[9]$ in the cycle notation.
   
 For an index $1\leq i \leq n$, we call it an excedance 
 of $\sigma$ 
if $\sigma(i)>i$, or a drop if $\sigma(i)<i$, or a fixed point if $\sigma(i)=i$.  Denote by $\exc(\sigma)$, $\drop(\sigma)$ and $\cyc(\sigma)$ the number of excedances, the number of drops and the number of cycles of $\sigma$, respectively. 
Let $D_n$ be the set of permutations without fixed points.  Then the derangement polynomials are defined by 
$$
d_n(x)=\sum_{\sigma \in D_n} x^{\exc(\sigma)},
$$
see \cite{Brenti-1990}. 

We can rely on the structure of  
forests of planted 0-1-2-plane trees
to give a combinatorial interpretations of the
$\gamma$-coefficients of the derangement polynomials,
and the
$q$-analogue with respect to the number of cycles, that is 
$$
d_n(x,y,q)=\sum_{\sigma \in D_n} x^{\exc(\sigma)}y^{\drop(\sigma)}q^{\cyc(\sigma)}.
$$
Notice that
a permutation without fixed points corresponds to a complete increasing binary tree without $\beta$-vertices whose left child is a leaf. A planted
increasing binary tree is said to be fully planted
if the root has a child that is not a leaf. 
By relabeling the root $1$ by $\beta$ 
and setting $\alpha=1$ in the $(a,b,\alpha,\beta)$-labeling, 
an increasing binary tree 
corresponding to a derangement can be decomposed into 
a forest of fully planted increasing binary trees 
for which the root of each tree is labeled by $\beta$. 
Then we take group action on a fully planted increasing 0-1-2-plane tree as follows. We label the root of each component
by $q$. If a non-root vertex has degree one,
then label it by $x+y$. A leaf is labeled by $xy$. Then the weight
of a forest $F$  of fully planted increasing 0-1-2-plane trees is
defined to be the product of all the grammatical labels of $F$, denoted
by $w(F)$. Then we get the following $\gamma$-expansion.
\begin{thm}
    For $n\geq 1$, we have
    \begin{equation}
        d_n(x,y,q)= \sum_{F}  w(F),
    \end{equation}
    where $F$ ranges over forests of fully planted increasing 0-1-2-plane trees
    on $[n]$. 
\end{thm}

\section{A labeling scheme for interior peaks}

In this section, we give two labeling schemes of increasing
binary trees in connection with interior peaks
of a permutation,  and we find combinatorial proofs of two
identities of Ji. 

Given a permutation $\sigma=\sigma_1\sigma_2\cdots\sigma_n$, 
 an index $i$ $(2\leq i \leq n-1)$ is called an interior peak if $\sigma_{i-1} < \sigma_{i}>\sigma_{i+1}$, and  we follow the notation  
$M(\sigma)$   in 
\cite{CF-2023-B} for  
the number of interior peaks of $\sigma$.

It turns out that the
number of interior peaks can be
read off from the $(a, b, \alpha, \beta)$-labeling
of the corresponding
increasing binary tree. 
More precisely,
a pair of sibling leaves 
labeled by $x$ and $y$ 
correspond to an interior peak of the
permutation. { {An $x$-leaf whose sibling is not a $y$-leaf   corresponds to an ascent  of the permutation. Likewise, a $y$-leaf
whose sibling is not
an $x$-leaf corresponds to
a descent of the permutation. }
Observe that the labels $a$ and $b$   play the role 
of preventing the first position
and the last position from being counted 
as interior peaks.

First, let us consider the $(\alpha,\beta)$-extension of
Stembridge's identity  \cite[Theorem 1.8]{Ji-2023-A}. 

\begin{thm}[Ji] \label{Ji1} For $n\geq 1$,
\begin{eqnarray}
   \lefteqn{ \sum_{\sigma\in S_n} (xy)^{M (\sigma)} \left(\frac{x+y}{2} \right)
    ^{n - 2{M}
    (\sigma)-1} \alpha^{\LR (\sigma) -1} \beta^{\RL(\sigma)-1}
    } \nonumber \\[6pt]
    & = & \sum_{\sigma\in S_n} x^{\des(\sigma)} y^{n-\des(\sigma)-1} 
    \left(\frac{\alpha+\beta}{2} \right)
    ^{\LR(\sigma)+\RL(\sigma)-2}.
    \label{ji-ab-s}
\end{eqnarray}
\end{thm}

The case for $n=1$ is trivial, so we assume that $n\geq 2$. To provide
a combinatorial interpretation of the above
relation, we shall give expansions of both sides
in terms of forests of  0-1-2-planted plane trees, and
will show that these two expansions are equinumerous,
that is, they amount to the same total weights.

To reformulate the above relation in terms of
trees, let $\mathcal{B}_n$ denote the
set of   increasing binary trees on $[n]$.
Given $T\in \mathcal{B}_n$ endowed with the $(a,b,\alpha,\beta)$-labeling, let $M(T)$ denote 
the number of  vertices of $T$ having two leaf
children. As used in \cite{Chen-Fu-2022},
$\xleaf(T)$ and
 $\yleaf(T)$ are referred to
  for the number of
$x$-leaves and the number of $y$-leaves 
of $T$.
Meanwhile, we write  $N_\alpha(T)$ and $N_\beta(T)$ for the
number of $\alpha$-vertices and the number of $\beta$-vertices of $T$, respectively. 

The relation (\ref{ji-ab-s}) 
can be split into two
parts. As for the left side, we have the 
following relation, where the
set of forests of planted 0-1-2-plane trees
on $[2,n]$ is denoted by $\mathcal{P}_n$.

\begin{thm} \label{TJ-1}
For $n\geq 1$, we have
    \begin{equation}\label{ji-ab-s-t-1}
    \sum_{T\in \mathcal{B}_n} (xy)^{M(T)} \left(\frac{x+y}{2} \right)
    ^{n - 2M(T)-1} \alpha^{N_\alpha (T) } \beta^{N_\beta(T) }
   =  \sum_{P \in \mathcal{P}_n } w(P) ,
    \end{equation}
where the sum ranges over the
set of forests of planted 0-1-2-plane 
trees on $[2,n]$ with the following labeling rules and 
$w(P)$ stands for the weight of $P$:
\begin{enumerate}
    \item A single root is labeled by
    $(x+y) \frac{\alpha + \beta}{2}$.

    \item The root of a component 
    with at least two  vertices
    is labeled by ${\alpha+\beta}$.
    
    \item A degree one vertex other than the root is labeled
    by $x+y$.
\end{enumerate}
\end{thm}

\noindent {\it Proof.} 
We begin with representing the sum on
the left side over permutations 
in terms of a sum over increasing binary trees. Let $T$ be an increasing binary tree of $\mathcal{B}_n$. 
 We say that a leaf is proper if it is neither 
an $a$-leaf nor a $b$-leaf. 
In view of the 
$(a,b,\alpha,\beta)$-labeling, 
$M(\sigma)$ corresponds to the
number of internal vertices having two proper
leaf children, whereas
$n-2M(\sigma)-1$ equals the number of internal vertices
having exactly one proper 
leaf child. 
Consequently, we are supposed 
to label $T$ by the following rules, which
we call the first modified $(a,b,
\alpha,\beta)$-labeling.

\begin{enumerate}
\item  
    
    Label the leftmost leaf by $a$ and
    label the rightmost leaf by $b$.

 \item 
    
    Any internal vertex on the path from the root to the $a$-leaf (other than the root) is labeled by $\alpha$. Any internal vertex on the path from the root to the $b$-leaf (other than the root) is labeled by $\beta$.

\item  For a pair of proper sibling leaves,
we label the left leaf by $x$ and the right leaf by $y$. 

        \item For a leaf whose sibling is
        not a proper leaf,
        we label it  by $(x+y)/2$, no matter whether it is on
        the left or on the right.

\end{enumerate}

For example, for the tree in Figure 
\ref{aabbf}, the first modified 
$(a,b,\alpha,\beta)$-labeling is demonstrated
in Figure \ref{aabbf2}.

\begin{figure}[!ht]
\begin{center}
\begin{tikzpicture}[scale=0.8]
\node [tn,label=90:$1$]{}[grow=down]
	[sibling distance=39mm,level distance=16mm]
    child {node [tn,label=180:{$4(\alpha)$}](four){}
       [sibling distance=26mm,level distance=13mm]
    child {node [tn,label=180:{$8(\alpha)$}](eight){}
     [sibling distance=14mm,level distance=13mm]
     child {node [tn1,label=-90:{}](eightl){}}
     child {node [tn1,label=0:{}](eightr){}}
     }
      child {node [tn,label=0:{$6$}](six){}
     [sibling distance=14mm,level distance=13mm]
     child {node [tn,label=180:{$9$}](nine){}
     child {node [tn1, label=-90:{}](ninel){}}
     child {node [tn1, label=-90:{}](niner){}}
     }
     child {node [tn1,label=0:{}](sixr){}}
     }
     }
     child {node [tn,label=0:{$2(\beta)$}](twoa){}
     [sibling distance=22mm,level distance=13mm]
      child {node [tn1,label=-90:{}](twoal){}}
    child {node [tn,label=0:{$3(\beta)$}](threea){}
       [sibling distance=22mm,level distance=13mm]
    child {node [tn,label=180:{$5$}](five){}
       [sibling distance=14mm,level distance=13mm]
      child {node [tn1,label=-90:{}](fivel){}}
    child {node [tn1,label=-90:{}](fiver){}}
    }
    child {node [tn,label=0:{$7(\beta)$}](sevena){}
         [sibling distance=14mm,level distance=13mm]
    child {node [tn1,label=-90:{}](sevenl){}}
    child {node [tn1,label=-90:{$b$}](sevenr){}}
    }
     }};
    \node [below=3pt] at (eightl){$a$}; 
     \node [below=-2pt] at (eightr){${x+y \over 2}$}; 
     \node [below=-2pt] at (sixr){${x+y \over 2}$};
     \node [below=3pt] at (ninel){$x$};
     \node [below=3pt] at (niner){$y$};
     \node [below=3pt] at (fivel){$x$}; 
     \node [below=3pt] at (fiver){$y$};
    \node [below=-2pt] at (sevenl){${x+y \over 2}$};
    \node [below=3pt] at (sevenr){$b$};
    \node [below=-2pt] at (twoal){${x+y \over 2}$};
    \node [tn1] at (eightr){};
    \node [tn1] at (sixr){};
    \node [tn1] at (sevenl){};
    \node [tn1] at (twoal){};
\end{tikzpicture}
\end{center}
\caption{The first modified 
$(a,b,\alpha,\beta)$-labeling. }
\label{aabbf2}
\end{figure}
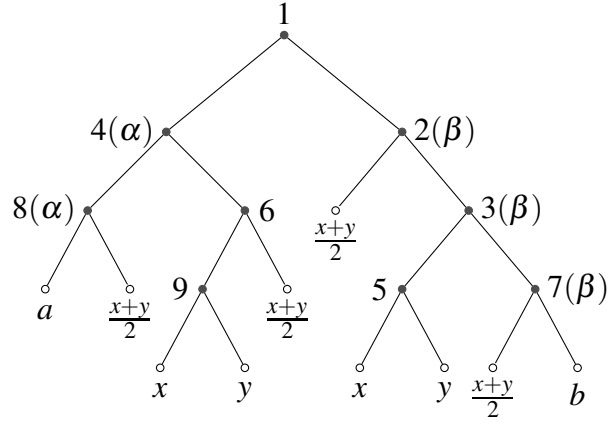

 We now process to compute the sum of
 weights over $\mathcal{B}_n$ 
 by utilizing
 supporting forests.
Let $F$ be a supporting forest, that is, a forest of planted increasing binary
trees on $[n]$. Let us characterize the 
set of trees $T$ in $\mathcal{B}_n$ 
with supporting forest $F$ on $[2,n]$. 
There are two choices for a planted increasing
binary tree in $F$ to belong to the
left side (with the root being
an $\alpha$-vertex) or the right side
(with the root being a $\beta$-vertex). 

For a single root, it may originate
from an $\alpha$-vertex in $T$ or a 
$\beta$-vertex in $T$. These two cases
lead to the sum of weights
\[ \frac{x+y}{2} \alpha  + \frac{x+y}{2}\beta =
    (x+y) \frac{\alpha+\beta}{2}. \]

For a component of $F$ containing at least 
two vertices, its root 
may originate from an $\alpha$-vertex
or a $\beta$-vertex, so the sum
of weights equals $\alpha+\beta$. 

Moreover, we can take a
  group action by swapping a proper leaf
with its sibling that is not a leaf. Keep 
in mind that the $a$-leaf and the $b$-leaf no 
longer appear in $F$. Let ${\rm orb}(F)$
denote the orbit of $F$ under this group action. Then let us compute the
sum of weights of $T$ with a supporting
forest in ${\rm orb}(F)$. This quantity
can be derived from a labeling of a 
representative of ${\rm orb}(F)$, 
that is a forest $P$ of planted 0-1-2-plane 
trees. 

Note that a proper left leaf with weight
$(x+y)/2$ is paired with a proper right
leaf with weight $(x+y)/2$, summing to
a weight $x+y$. The above considerations
suggest that we should comply with the
rules as stated
in the theorem. This completes the proof.
 \qed

Let us now turn to the sum on the right side
of (\ref{ji-ab-s}). A modification of the
$(a,b,\alpha, \beta)$-labeling is needed,
which we call the second modified
$(a,b,\alpha, \beta)$-labeling.
In this case, both 
the $\alpha$-vertices and the
$\beta$-vertices are labeled by $(\alpha+\beta)/2$. For example, 
Figure \ref{modifiedf} gives the modified labeling 
for the two trees in $\mathcal{B}_2$. 

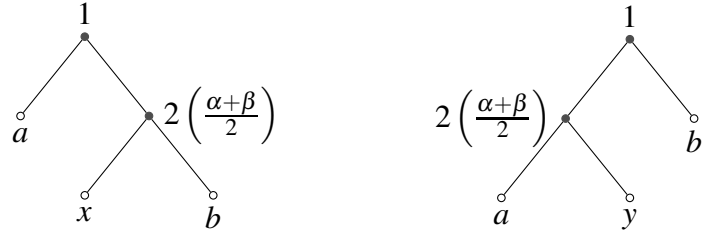
\begin{figure}[!ht]
\begin{center}
\begin{tikzpicture}[scale=0.9]
\node [tn,label=0:{}]{}[grow=down]
    child[grow=231] {node [tn1,label=-90:{$a$}](one){}}
     child[grow=309] {node [tn,label=0:{$2\left(\frac{\alpha+\beta}{2}\right)$}](){}
     child[grow=231] {node [tn1,label=-90:{$x$}](){}}
     child[grow=309] {node [tn1,label=-90:{$b$}](){}}};
     \node[above=2pt]{$1$};
    \end{tikzpicture} 
    \qquad \qquad \begin{tikzpicture}[scale=0.9]
\node [tn,label=180:{}]{}[grow=down]
    child[grow=231] {node [tn,label=180:{$ 2\left(\frac{\alpha+\beta}{2}\right)$}](one){}
              child[grow=231] {node [tn1,label=-90:{$a$}](){}}
     child[grow=309] {node [tn1,label=-90:{$y$}](){}}}
     child[grow=309] {node [tn1,label=-90:{$b$}](){}};
     \node[above=2pt]{$1$};
\end{tikzpicture}
\caption{The second modified 
$(a,b,\alpha, \beta)$-labeling.}
\label{modifiedf}
\end{center}
\end{figure}

At this point, the sum 
on the left of (\ref{ji-ab-s}) 
can be recast in terms
of the first modified $(a,b,\alpha,\beta)$-labeling for increasing
binary trees in $\mathcal{B}_n$, and we
are left with the task to establish
the following relation. 
 
\begin{thm}\label{thm-3.3}
For $n\geq 1$,
    \begin{equation}  \sum_{T\in \mathcal{B}_n} x^{\xleaf(T)} y^{\yleaf(T)} 
    \left(\frac{\alpha+\beta}{2} \right)
    ^{N_\alpha(T)+N_\beta(T)}=
    \sum_{P\in \mathcal{P}_n} w(P),
    \label{ji-ab-s-t-2}
\end{equation}
where the sum ranges over $\mathcal{P}_n$
as in Theorem \ref{TJ-1} and ditto the 
weight.
\end{thm}

\noindent {\it Proof.}
As before, we first consider the
supporting forest of a tree in $\mathcal{B}_n$, and consider which
trees in $\mathcal{B}_n $ share the
same supporting forest $F$. Let $T$ 
be an increasing binary tree in 
$\mathcal{B}_n $ with the supporting
forest $F$. 

For a single root in $F$, it may 
originate from an $\alpha$-vertex 
with a right leaf child labeled by $y$, 
or a $\beta$-vertex with a left leaf
child labeled by $x$. Given that all the
$\alpha$-vertices and $\beta$-vertices
are labeled by $(\alpha+\beta)/2$, 
the two cases contribute a total
weight of $(x+y)\frac{\alpha+\beta}{2}$,
in accordance with the labeling of $F$. 

For a component of $F$ containing at 
least two  vertices, its root 
may originate from an $\alpha$-vertex or
a $\beta$-vertex. Thus we get a total
weight of \[ \frac{\alpha+\beta}{2} +
\frac{\alpha+\beta}{2} = \alpha+\beta,\]
which coincides with the label of the
root of $F$.

For other leaves of $T$, we consider 
the group action that swaps a proper leaf
with its sibling that is an internal vertex. Strictly speaking, a proper $x$-leaf is paired with a proper $y$-leaf, giving
a total weight of $x+y$. This group action
gives rise to an orbit of $F$, which
can be represented by a forest of
planted 0-1-2-plane trees with weights
as designated in the theorem. This completes the proof. 
\qed

{{Next, let us recall the $(a,x,y,z)$-labeling
of an increasing binary tree  with the $x$-leaves, $y$-leaves and $z$-leaves
marking excedances, drops and fixed points of permutations respectively,}
see \cite{Chen-Fu-2024}:
\begin{enumerate}
    \item 
If a  $\beta$-vertex has a left leaf child, then this
child is labeled by $z$, signifying a fixed point. 
\item  The rest of the leaves are labeled 
in the same manner as the
$(a,b,\alpha,\beta)$-labeling with $a$  replaced by $x$ and 
$b$ replaced by $a$.
\end{enumerate} 
For example, 
with regard to the $(a,x,y,z)$-labeling, 
the increasing tree in Figure 1 corresponds to the following permutation in the cycle notation with the $(a,x,y,z)$-labels
attached: 
\begin{equation*}
(8\,y\,4\,x\,9\,y\,6 \, y \, 1 \, x) \;(2\,z )\; (5\, y \, 3 \, x) \; (7 \, z)\; a.
\end{equation*}

Clearly, a derangement corresponds to an increasing
binary tree without $z$-leaves. 

We finish with a combinatorial proof of the following identity due to Ji, where ${D}_n$ stands for the set of
derangements of $[n]$. 

 \begin{thm}[Ji]
 For $n\geq 1$,  
\begin{equation}\label{ji-i-3}
\sum_{\sigma \in  {S}_{n+1}}(-1)^{{\rm des}(\sigma)}{\left(\frac{1}{2}\right)}^{{\rm lrmin}(\sigma)+{\rm rlmin}(\sigma)-2}=\sum_{\sigma \in {  D}_{n}}(-1)^{{\rm exc}(\sigma)}.
\end{equation}
 \end{thm}

\noindent
{\it Proof. }
Let $T$ be a tree in $\mathcal{B}_{n+1}$ on $\{0,1,\ldots, n\}$ with the $(a,b,\alpha,\beta)$-labeling. 
Let $F$ be the supporting
forest of $T$. Consider the
set of trees that share the 
same supporting forest as $T$. First,  
we observe a cancellation property. Note that $F$ is a forest of planted
increasing binary trees on $[n]$.
On the other hand, 
we may regard $T$ as an
increasing binary tree endowed 
with the $(x,y)$-labeling for the
Eulerian polynomials, that is,
a left leaf is labeled by $x$ and a right
leaf is labeled by $y$. In the end,
we set  $x=-1$ and
$y=1$. 

We claim that a cancellation occurs when  
  $F$ contains a single root. If $F$ contains
  a single root, then $T$  has either an 
  $\alpha$-vertex with a $y$-leaf or 
    a $\beta$-vertex with an $x$-leaf. 
These two possibilities create a pair
of trees with the same supporting forest 
and opposite signs, here the sign of $T$ 
is determined by the parity of the number 
of $x$-leaves. Moreover, 
such a pair of trees possess the same quantity 
$$\LR(T)+\RL(T)-2,$$ 
and hence we are led to 
a cancellation in the
sum on the left of (\ref{ji-i-3}),
which implies that the sum 
can be reduced to 
$T$ whose supporting forests 
are fully planted. 

Note that the $(x,y)$-labels of $T$ are 
carried over to the forest $F$.
This means that
if two trees have the same supporting forest
(without single roots), then they must have the
same sign. Now we wish to  
compute the left side of (\ref{ji-i-3})
by classifying the supporting forests.
Clearly, a supporting forest of $k$ components
generates $2^k$ trees in $\mathcal{B}_{n+1}$. 

On the other hand, a supporting forest $F$
can be
viewed as an  increasing binary tree on $[n]$
by gluing the component together. Up to now, it remains to make use of the
fact that
the labels carried over are precisely the
same as the labels for the derangement 
polynomials with respect to the $(a,x,y,z)$-labeling,
except for the rightmost $y$-leaf.
 Thus  we may associate an  
  $x$-label with an excedance   and a $y$-label with a drop  of 
  the corresponding permutation. 
Finally, 
a special attention has to be paid to the 
rightmost $y$-leaf of $T$ in $\mathcal{B}_{n+1}$ 
subject to the $(a,x,y,z)$-labeling. Since $y$ is 
set to $1$ at last, there are no
worries. This completes the proof. 
\qed

To conclude, we remark that the above combinatorial argument yields a refinement 
of (\ref{ji-i-3}) by restricting the sum to 
\[ \LR(\sigma)+ \RL(\sigma)-2=k.\] Then the sum 
of the right side 
ranges over derangements with $k$ cycles.

\vskip 6mm \noindent{\large\bf Acknowledgment.}
We wish to thank the referee for 
helpful comments. This work was supported
by the National
Science Foundation of China.

\end{document}